\documentclass{article}
\usepackage[american]{babel}
\usepackage{amsfonts,amsmath, amsthm, amssymb}
\usepackage{xypic}
\xyoption{all}
\input{xypic}
\newdir{ >}{{}*!/-8pt/\dir{>}}
\def\R{\mathbb{R}}

\def\C{\mathbb{C}}

\def\Z{\mathbb{Z}}

\def\id{{\rm id}}\def\ker{{\rm ker}\,}\def\coker{{\rm coker}\,}
\def\ad{{\rm ad}}

\newtheorem{theo}{Theorem}
\newtheorem{defi}{Definition}
\newtheorem{Rem}{{\bf Remark}}
\newenvironment{rem}{\begin{Rem} \strut\\ \normalfont}{\end{Rem}}
\newtheorem{prop}{Proposition}

\newtheorem{cor}{Corollary}
\newtheorem{lem}{Lemma}
\begin{document}
\title{On Hopf $2$-algebras}

\author{Ya\"{e}l Fr\'egier\\
        Universit\'e du Luxembourg
\and  F. Wagemann\\
Universit\'e de Nantes}
\maketitle

\begin{abstract}
Our main goal in this paper is to translate the diagram below relating groups, 
Lie algebras and Hopf algebras to the corresponding $2$-objects, i.e. 
to categorify it.
This is done interpreting $2$-objects as crossed modules and showing the 
compatibility of the standard functors linking groups, Lie algebras 
and Hopf algebras with the concept of a crossed module. One outcome is the 
construction of an enveloping algebra of the string Lie algebra of Baez-Crans
\cite{BaeCra}, another is the clarification of the passage from crossed modules
of Hopf algebras to Hopf $2$-algebras.
\end{abstract}


\noindent keywords: crossed modules of Lie algebras, crossed modules of Hopf 
algebras, crossed modules of groups, $2$-groups, Lie $2$-algebras, Hopf 
$2$-algebras, string Lie algebra

\section*{Introduction}

One of the most impressive theorems in the theory of Lie groups is Lie's third 
theorem: the possibility to integrate a real or complex (finite-dimensional) 
Lie algebra in a unique way into a connected, simply connected Lie group. 
Algebraically, some aspects of this integration process are captured in the 
following diagram:

\hspace{4.5cm}
\xymatrix{
{\tt Lie} \ar[r]^{U} \ar[d] & {\tt ccHopf} \ar[d] \ar[l]^{P} \\ 
{\tt Grp} \ar[r]^{k[-]} \ar[u] & {\tt cHopf} \ar[u] \ar[l]^{\chi}}
\vspace{.5cm}

Here ${\tt Lie}$ is the category of Lie algebras over the field $k$ with
$k=\R$ or $k=\C$, ${\tt Grp}$ is the category of groups, supposed to be finite 
or connected algebraic (in which case we assume $k=\C$), 
when we apply the functor $k[-]$ of (regular) 
functions, ${\tt Hopf}$ is the category of $k$-Hopf algebras, and 
${\tt ccHopf}$ and ${\tt cHopf}$ are its subcategories of cocommutative resp.
commutative Hopf algebras. The functors $U$ and $P$ are those of the enveloping
algebra and of primitives, and $\chi$ is the functor of characters. The 
functors on the RHS of the diagram stipulate duality (linear duality in the 
case of finite dimensional Hopf algebras, 
restricted duality in the case of graded Hopf algebras with finite dimensional 
graded components,  or continuous duality in the case of complete topological 
Hopf algebras, depending on the context)
whereas those on the
LHS stipulate integration and derivation. Recall that the integration process 
consists in associating to a finite dimensional Lie algebra ${\mathfrak g}$
first its universal enveloping algebra $U{\mathfrak g}$. The latter is 
seen as the continuous dual of the completion of the space of functions on 
some formal group law (the space of ``point distributions''), see \cite{Ser}. 
Therefore in order to integrate
${\mathfrak g}$, one dualizes $U{\mathfrak g}$, which gives a completion
of the Hopf algebra of formal functions, 
and finally recovers the group as group-like elements or characters
in this Hopf algebra of functions. According to the properties of the 
initial Hopf algebra and the geometric requirements
of the reader, this gives then a formal group, an algebraic group or a finite 
group. For reference, the formal case may be found in \cite{Ser}.

On the other hand, in recent times, {\it categorification}, i.e. the passage
to categorified algebraic structures, plays a growing r\^ole in algebra and 
geometry. Here categorification means the replacement of the underlying sets
in some algebraic structure by categories, and maps between these sets by 
functors. For example, instead of regarding a group in ${\tt Sets}$, the 
category of sets (which is the ordinary concept of a group), one considers
a group object in the category ${\tt Cat}$ of (small) categories. In this way
one arrives at the notion of a $2$-group.  

There is recent intense mathematical research striving to understand how to
integrate Lie $2$-algebras into Lie $2$-groups, see for example \cite{BCSS},
\cite{Get}, \cite{Hen}, \cite{Woc}.
The purpose of this article is to establish the above diagram in the context
of $2$-groups and Lie $2$-algebras, which can be seen as some answer to the 
integration problem (see theorem $6$ in section $6$ in order to have a precise 
statement).
The diagram is obtained as the union of propositions
\ref{prop1} to \ref{prop4}. Thus our main result reads:

\begin{theo} \label{thm1}
The functors $U$, $P$, $k[-]$ and $\chi$, which we introduced above,
extend to the following diagram between categories of strict $2$-objects:

\vspace{.5cm}
\hspace{4.5cm}
\xymatrix{
{\tt 2-Lie} \ar[r]^{U} \ar[d] & {\tt 2-ccHopf} \ar[d] \ar[l]^{P} \\ 
{\tt 2-Grp} \ar[r]^{k[-]} \ar[u] & {\tt 2-cHopf} \ar[u] \ar[l]^{\chi}}
\vspace{.5cm}
\end{theo}

One aspect of this theorem is that it supplies the article \cite{FerLopNov}
with a huge amount of examples; namely, all strict Lie $2$-algebras give rise
to their kind of ${\rm cat}^1$-Hopf algebras (and, by the way, all semi-strict 
Lie $2$-algebras of \cite{BaeCra} can be strictified).  

The main method that we use in the proof of theorem \ref{thm1} is the 
translation of $2$-groups and Lie $2$-algebras
into a different algebraic structure, namely the structure of a 
{\it crossed module}. Let us explain this structure in the context of groups.
A {\it crossed module of groups} is a homomorphism of groups $\mu:M\to N$ 
together with an action $\alpha$ of $N$ on $M$ by automorphisms, denoted by 
$\alpha: m\mapsto ^nm$ for $n\in N$ and $m\in M$, such that 
\begin{itemize}
\item[(a)]  $\mu(^nm)\,=\,n\mu(m)n^{-1}$ and
\item[(b)]  $^{\mu(m)}m'\,=\,mm'm^{-1}$.  
\end{itemize}
It is well known, see \cite{Lod} \cite{MacLane} \cite{Por}, that the category
of strict $2$-groups and the category of
crossed modules of groups are equivalent. This is discussed in more 
detail in section $2$. Usually, this equivalence is formulated as an 
equivalence between $2$-categories, but in order to keep the abstract
categorical formalism to a minimum, we will stick to ordinary categories here.
We believe that the extension to $2$-categories is straight forward. 

Similarly, the category of strict Lie $2$-algebras and the category of
crossed modules of Lie algebras are equivalent (see \cite{BaeCra} and 
section $3$), and a similar statement also holds for cocommutative Hopf 
algebras \cite{FerLopNov}. 

Thus in order to achieve our task to transpose the above diagram into 
$2$-structures, we focus on the compatibility of the concept
of a crossed module with the above standard functors $U$, $P$, $k[-]$ and 
$\chi$. This is done in sections $4$ and $5$. In section $5$, due to the
duality coming into play when passing to the right bottom corner in the above
diagram, we explore the dual notion of crossed modules of Hopf algebra, namely,
{\it crossed comodules of Hopf algebras} (see definition \ref{def6}). 
We believe that this new algebraic structure gives an equivalent way of 
formulating Hopf $2$-algebras.

One outcome of the discussion of compatibility of the concept of a crossed 
module with standard functors is that the definition of a crossed module in 
\cite{FerLopNov} (taken over slightly generalized in our article as definition 
\ref{def1}) does not seem to be too far from the ``right'' categorification. 
Let us denote
such a crossed module of Hopf algebras by $\gamma:B\to H$.
Definition \ref{def1} 
imposes compatibility relations between the module structure (of $H$ on $B$)
and the Hopf algebra structure on $B$, namely, $B$ has to be an $H$-module 
algebra, an $H$-module coalgebra and the antipode of $B$ has to be a morphism
of $H$-modules. We show that these conditions are the natural reflection 
of the fact that
in the case of a crossed module of Lie algebras (see definition \ref{def5})
$\mu:{\mathfrak m}\to{\mathfrak n}$, the action has to be an action by 
derivations, i.e. for all $m, m'\in{\mathfrak m}$
and all $n\in{\mathfrak n}$: 
$$n\cdot[m,m']\,=\,[n\cdot m,m'] + [m,n\cdot m'],$$
and in the case of a crossed module of groups $\mu:M\to N$, the action has
to be an action by automorphisms of groups, i.e. for all $m,m'\in M$
and all $n\in N$:
$$^n(mm')\,=\,(^nm)(^nm').$$
Our main point is that this compatibility between the module structure 
(of $H$ on $B$) and the Hopf algebra structure on $B$ is necessary in 
case
one demands the crossed modules of Lie algebras and of groups to have similar
compatibility conditions. 

This is important to note. Indeed, the obstacle to define a crossed module of
associative algebras corresponding to a crossed module of Lie algebras by 
taking the functor $U$ term by term, is that the action is by derivations
of the associative product (see section $4$). Thus the term by term $U$-image 
of a crossed module of Lie algebras does not satisfy the 
condition of compatibility for a crossed module of associative 
algebras $\rho:R\to A$, which reads for all $a\in A$ and all $r,r'\in R$:
$$a(rr')\,=\,(ar)r',\,\,\,\,\,\,\,(rr')a\,=\,r(r'a),$$
see \cite{DIKL}: the naive belief that a crossed module of Hopf algebras is 
in particular a 
crossed module of associative algebras which is simultaneously a crossed module
of coassociative coalgebras is wrong.

On the other hand, the association of a crossed module of
associative algebras corresponding to a crossed module of Lie algebras may be
seen as a map on the level of cocycles whose induced map in cohomology
$$\phi:H^3({\mathfrak g},M^{\rm ad})\to HH^3(U{\mathfrak g},M)$$
is known to be an isomorphism. To our knowledge, no natural map on cocycles
inducing $\phi$ is known.   

In future work, we plan to drop the commutativity/cocommutativity condition 
in the definitions of a crossed (co)module of Hopf algebras
and to show in this more general framework the same equivalence of categories 
as in \cite{FerLopNov}.  

Let us note some by-products of our study: let $G$ be a connected, simply 
connected, complex simple Lie group and ${\mathfrak g}$ its Lie algebra. It is
well known that the de Rham cohomology group $H^3(G)$ is one dimensional and
isomorphic to $H^3({\mathfrak g})$, and that this space is generated by the 
Cartan cocycle $<[,],>$, which is manufactured from the Killing form $<,>$ and 
the bracket $[,]$ on ${\mathfrak g}$. In \cite{Wag}, the second named author 
gives an explicit 
crossed module $\mu:{\mathfrak m}\to{\mathfrak n}$ which represents the 
cohomology class of $<[,],>$ (via the bijection between equivalence classes of 
crossed modules $\mu:{\mathfrak m}\to{\mathfrak n}$ with fixed kernel 
$\ker(\mu)=V$ and fixed cokernel $\coker(\mu)={\mathfrak g}$ and 
$H^3({\mathfrak g},V)$). On the other
hand, the {\it string group} associated to $G$ is its $3$-connected cover, see
\cite{BCSS}, \cite{Hen1}. This string group is only defined up to homotopy and
cannot be realized as a strict, finite dimensional Lie $2$-group, but only as 
an infinite dimensional Lie $2$-group \cite{BCSS} or non-strict \cite{Woc2}. 
Applying the functor $U$ to our crossed 
module $\mu:{\mathfrak m}\to{\mathfrak n}$ gives by proposition \ref{prop1} 
a crossed module of Hopf algebras which is a natural candidate for an 
enveloping algebra of the string Lie algebra (see remark \ref{rem5}), 
explicitly:
\begin{cor} \label{cor1}
There is a crossed module of Hopf algebras 
$$\mu:S(N(0)^{\sharp})\to S(N(0)^{\sharp})\otimes_{\alpha}U{\mathfrak g}$$
which is a natural {\it algebraic} candidate for the enveloping algebra of the 
string Lie $2$-algebra in the sense that its underlying vector spaces are
at most countably infinite dimensional. 
\end{cor} 
This restriction on the dimension of the underlying vector spaces is important
when one wants to perform algebraic constructions with the enveloping algebra,
like for example pass to the dual crossed comodule.
On the other hand, it is well-known that a finite dimensional construction is 
not possible, cf \cite{Hoch}.
A similar statement holds for any finite dimensional semi-strict Lie 
$2$-algebra, see theorem $5$.
 
We leave open the question of defining an enveloping
algebra of a Lie $2$-algebra as a left adjoint $2$-functor (and refer to 
remark \ref{rem5} for the notations used in the preceding 
corollary).

A second by-product is a sort of Kostant's theorem for irreducible Hopf 
$2$-algebras in characteristic $0$ (see remark \ref{rem6}). 
Indeed, the ordinary Kostant theorem 
(asserting that an irreducible Hopf algebra in characteristic $0$ is isomorphic
to the universal enveloping algebra of its primitives) and
our techniques show how to find a crossed module of Lie algebras 
$\mu:{\mathfrak m}\to{\mathfrak n}$ for a given
crossed module of irreducible Hopf algebras $\gamma:B\to H$ such that 
$\gamma:B\to H$ is isomorphic to $U(\mu:{\mathfrak m}\to{\mathfrak n})$. 
\begin{cor}    \label{cor2}
An irreducible cocommutative Hopf $2$-algebra is equivalent to an enveloping
Hopf $2$-algebra, i.e. to a Hopf $2$-algebra of the form 
$U(\mu:{\mathfrak m}\to{\mathfrak n})$, where 
$\mu:{\mathfrak m}\to{\mathfrak n}$ is a crossed module of Lie algebras.
\end{cor}

A third by-product is a new approach to the equivalence between crossed modules
of Hopf algebras and Hopf $2$-algebras (called ${\rm cat}^1$-Hopf algebras in)
\cite{FerLopNov}. Namely, we use the functors $U$ and $P$ to reduce the 
problem to Lie algebras where the equivalence is well known. We get in this
way:

\begin{cor}  \label{cor3}
When restricting to the subcategory of irreducible cocommutative Hopf algebras,
the notions of crossed module and of ${\rm cat}^1$-Hopf algebra are equivalent.
\end{cor}

Open questions abound - let us state only two of them which will guide our
future research in this field: how to construct a cohomological interpretation
of crossed modules of Hopf algebras, inspired by the fact that crossed 
modules of Lie algebras, groups or associative algebras are classified (up
to equivalence) by $3$-cohomology classes ? How to quantize Hopf $2$-algebras,
or, in other words, how to deform our crossed modules of enveloping algebras
to get some quantum $2$-groups ? Some very partial answers to these questions 
are contained in remarks \ref{rem8} and \ref{rem9} in section $4$.\\   

{\bf Acknowledgements:} We thank the University of Luxembourg for financing 
stays of YF in Nantes and of FW in Luxembourg, where this research has been 
done. We thank the University of Technology Darmstadt for support during a 
small part of this work. We thank Chenchang Zhu for her interest in this work 
during the Seminar Sophus Lie in G\"ottingen in July 2009. We thank Urs 
Schreiber and Dimitry Roytenberg for answering a question concerning section 
$3$.

\section{Preliminaries}

Here we collect standard notations and definitions for studying Hopf algebras
coming from groups and Lie algebras, and their crossed modules. We introduce 
in definition \ref{def1} the notion of a 
 {\it crossed module of Hopf algebras} which is slightly more general than  
definition 12 given in \cite{FerLopNov} (see remark \ref{rem1}).

Let us denote by $\tau$ the twist in the symmetric monoidal category of vector
spaces over a commutative field $k$ of characteristic $0$. We refrain from 
formulating our entire paper relative to an arbitrary underlying 
symmetric monoidal category, but we believe that this generalization is 
straight forward. In geometric 
situations, we will always suppose $k=\C$.
All Lie and Hopf algebras are supposed to
be algebras over $k$. A Hopf algebra $H$ over $k$ is given by
$(H,\mu_H,\eta_H,\triangle_H,\epsilon_H,S_H)$, where $\mu_H$ is the associative
product, $\eta_H$ its unit, $\triangle_H$ the coassociative coproduct,
$\epsilon_H$ its counit, and $S_H$ is the antipode. In case the Hopf algebra 
we work with is clear from the context, we feel free to drop the corresponding
index (and write for example $\mu$ instead of $\mu_H$).

Let $H$ be a Hopf algebra. Then $H$ becomes a left $H$-module by the {\it 
adjoint action} which is defined for $h,k\in H$ as follows:
$$\ad_H(h)(k)\,=\,\mu\circ(\mu\otimes S)\circ(\id_H
\otimes\tau_{H,H})\circ(\triangle\otimes\id_H)(h\otimes k).$$

Recall the notions of an $H$-module algebra and an $H$-module coalgebra: 
given a Hopf algebra $H$ and an associative algebra $(A,\eta,\mu)$, $A$ is 
called a 
{\it left $H$-module algebra} in case the vector space underlying $A$ is a 
left $H$-module and the action is compatible with the algebra structure in 
the sense
that $\mu$ and $\eta$ are morphisms of $H$-modules (where $H$ acts on 
$A\otimes A$ and $k$ using $\triangle_H$ and $\epsilon_H$ respectively). 
Similarly, a coalgebra $(C,\epsilon,\triangle)$ is a {\it left $H$-comodule 
coalgebra} in case the vector space $C$ carries a left coaction of $H$ and 
$\epsilon$ and $\triangle$ are morphisms of $H$-comodules. In an analogous
way, one defines a {\it left $H$-comodule algebra} and a 
{\it left $H$-module coalgebra}. Some of the corresponding diagrams are spelt 
out in detail in the book of Kassel \cite{Kassel}.  
Everything extends in an obvious way to right modules. 

\begin{defi} \label{def1}
Let $\gamma:B\to H$ be a morphism of Hopf algebras. The morphism 
$\gamma:B\to H$ is a crossed module of Hopf algebras in case
\begin{itemize}
\item[(i)]  $(B,\phi_B)$ is a left $H$-module and a left $H$-module algebra 
and $H$-module coalgebra,
\item[(ii)] $\gamma\circ\phi_B\,=\,\mu_H\circ(\mu_H\otimes S_H)\circ(\id_H
\otimes\tau_{H,H})\circ(\triangle_H\otimes\gamma)$, or 
$\gamma\circ\phi_B\,=\,\ad_H\circ(\id_H\otimes\gamma)$ in shorthand notation, 
i.e. $\gamma$ is a 
morphism of $H$-modules where $H$ carries the $H$-module structure given by 
the adjoint action $\ad_H$,
\item[(iii)] $\phi_B\circ(\gamma\otimes\id_B)\,=\,\mu_B\circ(\mu_B\otimes S_B)
\circ(\id_B\otimes\tau_{B,B})\circ(\triangle_B\otimes\id_B)$, or 
$\phi_B\circ(\gamma\otimes\id_B)\,=\,\ad_B$ in shorthand 
notation (This is sometimes called 
{\it Peiffer identity}.)
\end{itemize}
\end{defi}

\begin{defi}
A morphism $(\rho,\sigma):(\gamma:B\to H)\to(\gamma':B'\to H')$
between crossed modules of Hopf algebras $\gamma:B\to H$ and $\gamma':B'\to H'$
is a pair of morphisms of Hopf algebras $\rho:B\to B'$ and $\sigma:H\to H'$
such that the following diagrams commute:

\vspace{.5cm}
\hspace{3cm}
\xymatrix{
B \ar[r]^{\rho} \ar[d]^{\gamma} & 
B' \ar[d]^{\gamma'} &  & B\otimes H \ar[r]^{\rho\otimes\sigma} \ar[d]^{\phi_B}
& B'\otimes H' \ar[d]^{\phi_B'} \\ 
H \ar[r]^{\sigma}  & H' & & B \ar[r]^{\rho} & B' }
\vspace{.5cm}
\end{defi}

It is easy to show that crossed modules of Hopf algebras form a category with
respect to this notion of morphisms.

\begin{rem} \label{rem1}
\begin{itemize}
\item[(a)] Observe that the definition is not autodual: $B$ carries only the
structure of an $H$-module - the definition does not demand a comodule
structure.
\item[(b)] Observe that one cannot associate naively a four term exact 
sequence of Hopf algebras to a crossed module of Hopf algebras: condition (ii) 
does not imply that the image of $\gamma$ is an associative ideal (while it 
is always a coideal).
\end{itemize}
\end{rem}

\begin{rem} \label{rem2}
 This definition differs from definition 12 given in \cite{FerLopNov}: we do 
not make any assumption on cocommutativity nor impose constraints coming from 
compatibility with cocommutativity (condition (i) of definition 12 in 
\cite{FerLopNov}). We do not ask the antipode $S_B$ to be a morphism of 
$H$-modules either (condition (iii) of definition 12 in \cite{FerLopNov}).
\end{rem}

Let us denote ${\tt Lie}$, ${\tt Hopf}$, ${\tt Grp}$, ${\tt alg Grp}$ the
categories of $k$-Lie algebras, resp. $k$-Hopf algebras, resp. groups, resp.
(connected) algebraic groups defined over $\C$. Speaking of algebraic groups 
will always imply that the ground field $k$ is $\C$, the field of complex 
numbers. Let us furthermore denote by
${\tt ccHopf}$ and ${\tt cHopf}$ the full subcategories of ${\tt Hopf}$
consisting of cocommutative resp. commutative Hopf algebras.
Recall the enveloping functor
$$U:{\tt Lie}\to{\tt Hopf},$$
which associates to a Lie algebra ${\mathfrak g}$ its universal enveloping
algebra $U{\mathfrak g}$; the functor of primitives
$$P:{\tt Hopf}\to{\tt Lie},$$
which associates to a Hopf algebra its Lie algebra of primitive elements
(i.e. the $h\in H$ such that $\triangle h=1\otimes h+h\otimes 1$); the
functor of regular functions
$$k[-]:{\tt alg Grp}\to{\tt cHopf},$$
which associates to a connected algebraic group $G$ its Hopf algebra of regular
functions $k[G]$; the functor of characters
$$\chi:{\tt Hopf}\to{\tt Grp},$$
which associates to a Hopf algebra $H$ the group of its characters
$\chi(H)$ (i.e. of algebra morphisms from $H$ to $k$).
Recall that a commutative Hopf algebra $H$ which is finitely generated as an
algebra gives rise to an affine algebraic group, see \cite{Swe}, section 6.3, 
p. 123. The associated algebraic group is $\chi(H)$. 
Recall furthermore that the functors $U$ and $P$ are equivalences of categories
in characteristic $0$ when one restricts to the full 
subcategory of {\it irreducible} cocommutative Hopf algebras (see \cite{Swe} 
theorem 13.0.1, p. 274). 

\section{On $2$-groups}

The subject of this section is the categorification of the notion of a group, 
a categorified group being a $2$-group. This matter is well known and we refer 
to \cite{MacLane}, chapter XII, \cite{Lod} and \cite{Por}.

To categorify an algebraic notion, one looks at this kind of object not in the
category of sets, but in the category of (small) categories. For example, a
categorified group, or $2$-group, is a group object in the category of
categories. Amazingly, this is the same as a category object in the category of
groups, see \cite{MacLane} p. 269. 

\begin{defi} \label{def2}
A $2$-group is a category object in the category ${\tt Grp}$, i.e. it is 
the data of two groups $G_0$, the group of objects, and
$G_1$, the group of arrows, together with group homomorphisms $s,t:G_1\to G_0$,
source and target, $i:G_0\to G_1$, inclusion of identities, and
$m:G_1\times_{G_0}G_1\to G_1$, the categorical composition (of arrows) 
which satisfy the usual axioms of a category.
\end{defi}

We should emphasize, however, that the $2$-groups introduced here
are {\it strict} $2$-groups, i.e. in the categorified version all laws are 
verified as equations. In other categorifications, one relaxes some or all 
laws to hold only
up to natural transformation, transformations which should then satisfy
coherence conditions. This leads to much more general notions (like coherent 
$2$-groups, weak $2$-groups, see \cite{BaeLau}), but does not occupy us here.

Recall the notion of a morphism between strict $2$-groups:

\begin{defi}
A morphism $F:(G_0,G_1)\to(H_0,H_1)$ between $2$-groups $(G_0,G_1)$ and 
$(H_0,H_1)$ is a functor internal to the category ${\tt Grp}$, i.e. the data
of morphisms of groups $F_0:G_0\to H_0$ and $F_1:G_1\to H_1$ such that the 
following diagrams commute:

\hspace{2cm}
\xymatrix{
G_1\times_{G_0}G_1 \ar[r]^{F_1\times F_1} \ar[d]^{m_G} & 
H_1\times_{H_0}H_1 \ar[d]^{m_H} &  & G_1 \ar[r]^{s,t} \ar[d]^{F_1}
& G_0 \ar[r]^{i} \ar[d]^{F_0} & G_1 \ar[d]^{F_1} \\ 
G_1 \ar[r]^{F_1}  & H_1 & & H_1 \ar[r]^{s,t} & H_0 \ar[r]^{i} & H_1 }
\vspace{.5cm}
\end{defi}
   
With this notion of morphisms, strict $2$-groups form a category. $2$-groups 
also form a $2$-category (cf \cite{Por}, \cite{MacLane}), but we will stick to 
the easiest categorical framework in which our article works. 

The notion of a crossed module of groups has already been defined in the 
introduction. We recall it here for the convenience of the reader:

\begin{defi}
A crossed module of groups is a homomorphism of groups $\mu:M\to N$ 
together with an action $\alpha$ of $N$ on $M$ by automorphisms, denoted by 
$\alpha: m\mapsto ^nm$ for $n\in N$ and $m\in M$, such that 
\begin{itemize}
\item[(a)]  $\mu(^nm)\,=\,n\mu(m)n^{-1}$ and
\item[(b)]  $^{\mu(m)}m'\,=\,mm'm^{-1}$.  
\end{itemize} 
\end{defi}

Recall the notion of a morphism of crossed modules of groups:

\begin{defi}
A morphism $(\rho,\sigma):(\mu:M\to N)\to(\mu':M'\to N')$ between crossed 
modules of groups $(\mu:M\to N)$ and $(\mu':M'\to N')$ is a pair of group 
homomorphisms $\rho:M\to M'$ and $\sigma:N\to N'$ such that the following 
diagrams commute:

\hspace{3cm}
\xymatrix{
M \ar[r]^{\rho} \ar[d]^{\mu} & 
M' \ar[d]^{\mu'} &  & M\times N \ar[r]^{\rho\times\sigma} \ar[d]^{\alpha}
& M'\times N' \ar[d]^{\alpha'} \\ 
N \ar[r]^{\sigma}  & N' & & M \ar[r]^{\rho} & M' }
\vspace{.5cm}
\end{defi}

Crossed modules form a category with this notion of morphisms. They also form
a $2$-category, as explained for example in \cite{Por}. 
The following theorem can be found in \cite{MacLane}, chapter 
XII, \cite{Lod}, and in much more detail in \cite{Por}.

\begin{theo} \label{thm2}
The categories of $2$-groups and of crossed modules of groups are equivalent.
\end{theo}


\begin{rem}\label{rem4}
One interesting point about this theorem is that the categorical 
composition $m:G_1\times_{G_0}G_1\to G_1$ does not represent additional 
structure, but is already encoded in the group law of $G_1$, namely, one has
$$g\circ f\,:=\,m(f,g)\,=\,f(i(b))^{-1}f,$$
where $t(f)=b$; this formula (which involves only the group multiplication in 
$G_1$ on the RHS) is shown, for example, in chapter 
XII of \cite{MacLane}. Thus the data of two groups $G_0$, $G_1$ and morphisms
$s,t:G_1\to G_0$ and $i:G_0\to G_1$ satisfying the usual axioms of source, 
target and object inclusion in a category is already equivalent to the data of 
a crossed module. 
\end{rem}

\section{On Lie $2$-algebras}

Here we recall in the same way as in the previous section the correspondence
between crossed modules of Lie algebras and strict Lie $2$-algebras. Note that 
this correspondence does not occur as stated in the literature, but is treated
for semi-strict Lie $2$-algebras in \cite{BaeCra}. 

In order to define a Lie $2$-algebra, we first define a $2$-vector 
space (over the ground field $k$).

\begin{defi} \label{def3}
A $2$-vector space is a category object in the category ${\tt Vect}$ of vector 
spaces, i.e. it is the data of a vector space $V_1$ of arrows, a vector space
$V_0$ of objects, and of linear maps $s,t:V_1\to V_0$, $i:V_0\to V_1$ and
$m:V_1\times_{V_0}V_1\to V_1$ satisfying the usual axioms of a category. 
\end{defi}

$2$-vector spaces come together with their natural notion of morphisms:

\begin{defi}
A morphism $F:V\to V'$ between $2$-vector spaces $V$ and $V'$ is a linear
functor, i.e. the data of two linear maps $F_0:V_0\to V_0'$ and 
$F_1:V_1\to V_1'$ such that the following diagrams commute:

\vspace{.5cm}
\hspace{2cm}
\xymatrix{
V_1\times_{V_0}V_1 \ar[r]^{F_1\times F_1} \ar[d]^{m_V} & 
V_1'\times_{V_0'}V_1' \ar[d]^{m_{V'}} &  & V_1 \ar[r]^{s,t} \ar[d]^{F_1}
& V_0 \ar[r]^{i} \ar[d]^{F_0} & V_1 \ar[d]^{F_1} \\ 
V_1 \ar[r]^{F_1}  & V_1' & & V_1' \ar[r]^{s,t} & V_0' \ar[r]^{i} & V_1' }
\vspace{.5cm}
\end{defi}

Once again, the fact that we are dealing with strict $2$-vector spaces
implies that we have a category of $2$-vector spaces (with the above notion of 
morphisms).    

It is shown in lemma 6 of \cite{BaeCra} that once again, the categorical 
composition $m$ is redundant data, i.e. can be recovered by the vector space 
structure. This can also be seen as a special case of the corresponding result 
of the previous section.
More explicitly, writing elements $f,g\in V_1$ as $f=\vec{f}+i(s(f))$ and 
$g=\vec{g}+i(s(g))$, one has 
$$g\circ f\,:=\,m(f,g)\,=\,i(s(f))+\vec{f}+\vec{g}$$
(Note that we have the usual convention for the composition of maps, while 
Baez and Crans have the categorical convention.)

\begin{defi} \label{def4}
A Lie $2$-algebra is a category object in the category ${\tt Lie}$ of Lie 
algebras, i.e. it is the data of a Lie algebra ${\mathfrak g}_1$ of arrows, 
a Lie algebra ${\mathfrak g}_0$ of objects, and of Lie algebra morphisms 
$s,t:{\mathfrak g}_1\to{\mathfrak g}_0$, 
$i:{\mathfrak g}_0\to{\mathfrak g}_1$ and
$m:{\mathfrak g}_1\times_{{\mathfrak g}_0}{\mathfrak g}_1\to {\mathfrak g}_1$ 
satisfying the usual axioms of a category. 
\end{defi}

Lie $2$-algebras come together with the appropriate notion of morphisms:

\begin{defi}
A morphism $F:{\mathfrak g}\to{\mathfrak h}$ between (strict) Lie $2$-algebras
${\mathfrak g}$ and ${\mathfrak h}$ is a morphism of the underlying $2$-vector 
spaces such that the linear maps $F_0$ and $F_1$ are morphisms of Lie algebras.
\end{defi}

Comparing this definition to the corresponding one in \cite{BaeCra} (def. 23), 
one sees that $F_2$ is always the identity in our setting.  

In particular, a Lie $2$-algebra is a $2$-vector space and the categorical 
composition 
$m:{\mathfrak g}_1\times_{{\mathfrak g}_0}{\mathfrak g}_1\to {\mathfrak g}_1$
can be recovered from the $2$-vector space structure (the fact that the 
so-defined $m$ is a morphism of Lie algebras is easily verified).  

Let us also recall the definition of a crossed module of Lie algebras, for more
information on these see \cite{Wag}. 

\begin{defi}  \label{def5}
A crossed module of Lie algebras is a homomorphism of Lie algebras
$\mu:{\mathfrak m}\to{\mathfrak n}$ together with an action, denoted 
$\alpha(m,n)\,=\,n\cdot m$, of ${\mathfrak n}$
on ${\mathfrak m}$ by derivations such that for all $m,m'\in{\mathfrak m}$ and 
all $n\in{\mathfrak n}$
\begin{itemize}
\item[(a)]  $\mu(n\cdot m)\,=\,[n,\mu(m)]$ and
\item[(b)]  $\mu(m)\cdot m'\,=\,[m,m']$.
\end{itemize}
\end{defi}

\begin{defi}
A morphism $(\rho,\sigma):(\mu:{\mathfrak m}\to{\mathfrak n})\to
(\mu':{\mathfrak m}'\to{\mathfrak n}')$ of crossed morphisms of Lie algebras 
is a pair of morphisms of Lie algebras $\rho:{\mathfrak m}\to{\mathfrak m}'$ 
and $\sigma:{\mathfrak n}\to{\mathfrak n}'$ such that the following diagrams 
commute: 

\hspace{3cm}
\xymatrix{
{\mathfrak m} \ar[r]^{\rho} \ar[d]^{\mu} & 
{\mathfrak m}' \ar[d]^{\mu'} &  & {\mathfrak m}\times{\mathfrak n}
\ar[r]^{\rho\times\sigma} \ar[d]^{\alpha}
& {\mathfrak m}'\times{\mathfrak n}' \ar[d]^{\alpha'} \\ 
{\mathfrak n} \ar[r]^{\sigma}  & {\mathfrak n}' & & {\mathfrak m} 
\ar[r]^{\rho} & {\mathfrak m}' }
\vspace{.5cm}
\end{defi}

\begin{theo}  \label{thm3}
The categories of strict Lie $2$-algebras and of crossed modules of Lie 
algebras are equivalent.
\end{theo}

\begin{proof} Given a Lie $2$-algebra $s,t:{\mathfrak g}_1\to{\mathfrak g}_0$,
the corresponding crossed module is defined by
$$\mu:=t|_{\ker(s)}\,:\,{\mathfrak m}:=\ker(s)\to{\mathfrak n}:=
{\mathfrak n}.$$
The action of ${\mathfrak n}$ on ${\mathfrak m}$ is given by
$$n\cdot m\,:=\,[i(n),m],$$
for $n\in{\mathfrak n}$ and $m\in{\mathfrak m}$ (where the bracket is taken in 
${\mathfrak g}_1$). This is well defined and an action by derivations.
Axiom (a) follows from 
$$\mu(n\cdot m)=\mu([i(n),m])=[\mu\circ i(n),\mu(m)]=[n,\mu(m)].$$
Axiom (b) follows from
$$\mu(m)\cdot m'=[i\circ\mu(m),m']=[i\circ t(m),m']$$
by writing $i\circ t(m)=m+r$ for $r\in\ker(t)$ and by using that $\ker(t)$ and 
$\ker(s)$ in a Lie $2$-algebra commute (shown in Lemma \ref{lem1} after the 
proof).

On the other hand, given a crossed module of Lie algebras 
$\mu:{\mathfrak m}\to{\mathfrak n}$, associate to it
$$s,t\,:\,{\mathfrak n}\ltimes{\mathfrak m}\to{\mathfrak n}$$
by $s(n,m)=n$, $t(n,m)=\mu(m)+n$, $i(n)=(n,0)$, where the semi-direct product
Lie algebra ${\mathfrak n}\ltimes{\mathfrak m}$ is built from the given action 
of ${\mathfrak n}$ on ${\mathfrak m}$ by $[(n_1,m_1),(n_2,m_2)]:=(
[n_1,n_2]_{\mathfrak n},n_1\cdot m_2-n_2\cdot m_1)$.
Let us emphasize that this semi-direct 
product concerns {\it only} the Lie algebra structure on ${\mathfrak n}$ and 
the ${\mathfrak n}$-module structure on ${\mathfrak m}$; the bracket on 
${\mathfrak m}$ is lost. (Nevertheless, the bracket on ${\mathfrak m}$ in a 
crossed module $\mu:{\mathfrak m}\to{\mathfrak n}$ is encoded in the action and
the morphism $\mu$ by axiom (b).) The composition of arrows is already
encoded in the underlying structure of $2$-vector space, as remarked before 
the statement of the theorem.
\end{proof}

\begin{lem}  \label{lem1}
$[\ker(s),\ker(t)]\,=\,0$ in a strict Lie $2$-algebra.
\end{lem}    

\begin{proof} 
The fact that the composition of arrows is a homomorphism of Lie algebras
gives the following ``middle four exchange'' property
$$[g_1,g_2]\circ[f_1,f_2]\,=\,[g_1\circ f_1,g_2\circ f_2]$$
for composable arrows $f_1,f_2,g_1,g_2\in{\mathfrak g}_1$. 
Now suppose that $g_1\in\ker(s)$ and $f_2\in\ker(t)$. Then denote by $f_1$ and 
by $g_2$ the
identity (w.r.t. the composition) in $0\in{\mathfrak g}_0$. As these are 
identities, we have $g_1=g_1\circ f_1$ and $f_2=g_2\circ f_2$. On the other 
hand, $i$ is a morphism of Lie algebras and sends $0\in{\mathfrak g}_0$ to
the $0\in{\mathfrak g}_1$. Therefore we may conclude
$$[g_1,f_2]\,=\,[g_1\circ f_1,g_2\circ f_2]\,=\,[g_1,g_2]\circ[f_1,f_2]\,=\,
0.$$
\end{proof}

\section{Crossed modules of Lie and Hopf algebras}

The goal of this section is to study the compatibility of the notion of crossed
module with the standard functors $U$ and $P$ between the categories of
Lie algebras and cocommutative Hopf algebras. The main proposition shows how to
associate a Hopf $2$-algebra to a Lie $2$-algebra.

\begin{lem} \label{lem2}
An action of a Lie algebra ${\mathfrak g}$ by derivations (of the Lie bracket)
on a Lie algebra ${\mathfrak h}$ extends to an action by derivations (of the
associative product) of $U{\mathfrak g}$ on $U{\mathfrak h}$.
\end{lem}

\begin{proof} This is prop. 2.4.9 (i), p. 81, in \cite{Dix}.\end{proof}

\begin{lem}  \label{lem3}
Let ${\mathfrak g}$ be a Lie algebra. An associative algebra $A$ is an 
$U{\mathfrak g}$-module algebra if and only if the vector space $A$ is a 
${\mathfrak g}$-module such that elements of ${\mathfrak g}$ act by derivations
(of the associative product of $A$).
\end{lem}

\begin{proof} This is lemma V.6.3, p. 108, in \cite{Kassel}. \end{proof}

\begin{prop} \label{prop1}
The functor $U$ (see section $1$) extends to a functor (still denoted $U$) 
from the category of
crossed modules of Lie algebras to the catoegory of crossed modules of
(cocommutative) Hopf algebras.
\end{prop}

\begin{proof} 
Let $\mu:{\mathfrak m}\to{\mathfrak n}$ be a crossed module of Lie 
algebras. By functoriality, we get a homomorphism of associative algebras
$\gamma:=U\mu:U{\mathfrak m}\to U{\mathfrak n}$. By the previous two lemmas,
$U{\mathfrak n}$ acts (by derivations of the associative product)
on $U{\mathfrak m}$ and with this action, $U{\mathfrak m}$ becomes an 
$U{\mathfrak n}$-module algebra. Let us show that $U{\mathfrak m}$ is also an 
$U{\mathfrak n}$-module coalgebra, i.e. that the coproduct 
$\triangle_{U{\mathfrak m}}$ and the counit $\epsilon_{U{\mathfrak m}}$ are 
$U{\mathfrak n}$-module homomorphisms. Recall for this the 
action of $n\in U{\mathfrak n}$ on the tensor product
$$n\cdot m\otimes m'\,:=\,
\triangle_{U{\mathfrak n}}n\cdot m\otimes m',$$
where $m,m'\in U{\mathfrak m}$. Now for primitive elements $m\in {\mathfrak m}
\subset U{\mathfrak m}$ and $n\in {\mathfrak n}\subset U{\mathfrak n}$, one
has
$$\triangle_{U{\mathfrak n}}n\cdot\triangle_{U{\mathfrak m}}m\,=\,
\triangle_{U{\mathfrak m}}(n\cdot m),$$
because ${\mathfrak n}$ acts trivially on $k$. The general case is obtained
using induction and the fact that the coproducts are algebra morphisms.
The counit $\epsilon:U{\mathfrak m}\to k$ is clearly a morphism.

It remains to show properties (ii) and (iii). These two follow from the 
properties (a) and (b) of a crossed module (see definition \ref{def5}). 
Indeed, on 
primitives, identity (ii) {\it is} identity (a) and identity (iii) {\it is}
identity (b). The general case follows from induction. 

One word about the morphism-side of the statement: a morphism of crossed 
modules $(\rho,\sigma):(\mu:{\mathfrak m}\to{\mathfrak n})\to
(\mu':{\mathfrak m}'\to{\mathfrak n}')$ is sent to a morphism of crossed 
modules of Hopf algebras $(U\rho,U\sigma):(U\mu:U{\mathfrak m}\to
U{\mathfrak n})\to(U\mu':U{\mathfrak m}'\to U{\mathfrak n}')$. This renders 
the first diagram in Definition $2$ commutative. In order to render the second
one commutative, one should use the natural isomorphism 
$\zeta_{{\mathfrak m},{\mathfrak n}}:U({\mathfrak m}\times{\mathfrak n})\to
U{\mathfrak m}\otimes U{\mathfrak n}$ and the actions 
$\phi_{U{\mathfrak m}}:U{\mathfrak n}\otimes U{\mathfrak m}\to U{\mathfrak m}$
given by $\phi_{U{\mathfrak m}}=\alpha\circ
\zeta_{{\mathfrak m},{\mathfrak n}}^{-1}$, where $\alpha$ is the action of
${\mathfrak n}$ on ${\mathfrak m}$ from the crossed module.
\end{proof}

\begin{prop}  \label{prop2}
The functor $P$ (see section $1$) extends to a functor (still denoted $P$)
from the category of crossed modules of Hopf algebras to the category of
crossed modules of Lie algebras.
\end{prop}

\begin{proof} 
Let $\gamma:B\to H$ be a crossed module of Hopf algebras (definition
\ref{def1}). The set of primitives $P(B)$ and $P(H)$ of $B$ and $H$ are Lie 
algebras,
and the restriction of $\gamma$ to $P(B)$ is a Lie algebra morphism
$\gamma:P(B)\to P(H)$. By hypothesis, we have a map $\zeta:H\otimes B\to B$
which is an action of $H$ on $B$.
Restrict it to $\bar{\zeta}:H\otimes P(B)\to
B$. $\bar{\zeta}$ takes its values in $P(B)$, because $(\star)$ the 
coproduct of $B$ is a morphism of $H$-modules:
\begin{eqnarray*}
\Delta (h\cdot b) & \overset{(\star)}{=} & \Delta(h)\cdot \Delta b\\
                 & = & (\sum h'\otimes h'')\cdot(\sum b'\otimes b'')\\
                 & \overset{b\in P(B)}{=} & \sum(h'\cdot 1\otimes h''\cdot b +
h'\cdot b\otimes h''\cdot 1)\\
                 & = & \sum(\epsilon(h')\cdot 1\otimes h''\cdot b +
h'\cdot b\otimes \epsilon(h'')\cdot 1)\\
                 & = & 1\otimes (\sum(\epsilon(h')h'')\cdot b) +
((\sum\epsilon(h'')h')\cdot b)\otimes 1)\\
                 & = & 1\otimes (h\cdot b) + (h\cdot b)\otimes 1)
\end{eqnarray*} 
Then we may restrict $\bar{\zeta}$ further to 
$\tilde{\zeta}:P(H)\otimes P(B)\to P(B)$.
The ``associativity'' of $\tilde{\zeta}$ is clear (and implies the property of 
a Lie algebra action), but we have to show that
the action $\tilde{\zeta}$ is by derivations. This follows from the hypothesis
that the multiplication of $B$ is a morphism of $H$-modules.

As we have already mentioned, properties (ii) and (iii) of a crossed module 
of Hopf algebras imply, when restricted to the primitives, properties (a) and 
(b) of a crossed module of Lie algebras.

The morphism-side of the statement is straight forward.\end{proof}

\begin{rem}  \label{rem6}
Let us deduce from the above propositions some equivalence between irreducible 
cocommutative Hopf $2$-algebras over a field of characteristic $0$ and 
enveloping $2$-algebras.

Indeed, it is well known that the functors $U$ and $P$ are equivalences of 
categories in characteristic $0$ when one restricts the category of Hopf 
algebras to the full subcategory of irreducible cocommutative Hopf 
algebras (see \cite{Swe} theorem 13.0.1, p. 274). 
Call a Hopf $2$-algebra {\it irreducible} in case 
its corresponding crossed module of Hopf algebras is composed of irreducible 
Hopf algebras. One deduces then from proposition \ref{prop1}:

\begin{theo}  \label{theo5}
An irreducible cocommutative Hopf $2$-algebra is isomorphic to an enveloping
Hopf $2$-algebra, i.e. to a Hopf $2$-algebra of the form 
$U(\mu:{\mathfrak m}\to{\mathfrak n})$, where 
$\mu:{\mathfrak m}\to{\mathfrak n}$ is a crossed module of Lie algebras.
\end{theo}
\end{rem}

\begin{rem}  \label{rem7}
In a very similar manner as the proofs of propositions \ref{prop1} and 
\ref{prop2}, 
it is easy to show that in case $\mu:M\to N$ is 
a crossed module of groups, the induced morphism between group
algebras $\mu:kM\to kN$ is a crossed module of Hopf algebras. The fact that
$kM$ is a $kN$-module algebra comes from the fact that the $N$-action on $M$
is by automorphisms, and that $kM$ is a $kN$-module coalgebra is always the 
case when one linearizes an action of a group on a set, see \cite{Kassel} 
p.203.
Taking group-like elements is the way back to the crossed module of groups.

For infinite groups, the corresponding group algebras do not (in general)
reflect all features of the given group. Therefore one often passes to 
topological versions, like the $C^*$-algebra associated to a group. The
correspondence between crossed modules of groups and crossed modules of their
group algebras respects this kind of topological versions. In case one sees
the $C^*$-algebra associated to a group not as some completion of the group 
algebra, but as a space of characters on the group, it is closer to the 
constructions in section $5$, but the same remark applies. 

As an example, consider the group $S^1$ and the dense subgroup $\Z$ which is 
the image of the embedding $\Z\to S^1$, $k\mapsto e^{i\lambda k}$ where 
$\lambda\,/\,2\pi$ is an irrational number. 
The quotient $S^1\,/\,\Z$ is an abelian
group, but carries the discrete topology, thus the associated $C^*$-algebra is
trivial. It is explained in \cite{BloWei} (see also references therein) that
a way to associate meaningful $C^*$-algebras to this example goes under the
name {\it Hopfish algebras}. We will not go into the definition of a Hopfish 
algebra and refer to {\it loc. cit.}. 

The point we want to make here (we owe this remark to Chenchang Zhu) is that
another way around the problem of associating a $C^*$-algebra to the quotient
$S^1\,/\,\Z$ is to take the crossed module of $C^*$-algebras associated to
the crossed module of groups $\Z\to S^1$. Its corresponding $2$-group (see 
section $1$) is then $S^1\ltimes\Z\to S^1$ (which may be regarded as a 
groupoid). In other words, this crossed module of $C^*$-algebras
replaces the $C^*$-algebra of the groupoid 
$S^1\ltimes\Z\to S^1$. A more precise link between the $C^*$-algebra of the 
groupoid and the crossed module of $C^*$-algebras needs to be investigated.
Note that this framework applies for any group $G$ and 
for any normal subgroup $N$ to the crossed module $N\to G$ or the corresponding
$2$-group/groupoid $G\ltimes N\to G$.   
\end{rem}

\begin{rem}  \label{rem8}
Let us comment on the possibility of a ``standard cohomological framework'' for
crossed modules of Hopf algebras. In the case of Lie algebras, (discrete) 
groups or associative algebras, taking sections one may associate to a crossed 
module a $3$-cohomology class which expresses the equivalence class of the 
given crossed module. As we stated before, this scheme needs modification in
the context of Hopf algebras. One way around this problem is to transfer using
the functor $U$ and proposition \ref{prop1} the equivalence relation for 
crossed 
modules of Lie algebras and their cohomology classes directly to Hopf algebras.
This may not be the most canonical way to do so, but it sets at least 
compatibility conditions which one might want to impose on some cohomological 
description of crossed modules of Hopf algebras. 
\end{rem} 

\begin{rem} \label{rem9}
Let us emphasize that Baez and Crans \cite{BaeCra}, section $6$, thought of 
$\{{\mathfrak g}_{\hbar}\}_{\hbar}$ as 
a one-parameter deformation of Lie $2$-algebras of the ``trivial'' 
Lie $2$-algebra ${\mathfrak g}$. From this point of view, the family of
crossed modules of Hopf algebras $U{\mathfrak g}_{\hbar}$ may be regarded as a
quantum $2$-group, i.e. a non-trivial deformation of Hopf $2$-algebras. 

In order to explain in which sense this deformation
is non-trivial, we need the $2$-categorical structure. Each semi-strict Lie 
$2$-algebra ${\mathfrak g}_{\hbar}$ is non-strict and only (non-trivially) 
equivalent to a strict one (``strictification''),
while ${\mathfrak g}_{0}$ is strict. 
\end{rem}
    
\section{Crossed (co)modules of groups and Hopf algebras}

The goal of this section is to study the compatibility of the notion of crossed
module with the standard functors $k[-]$ and $\chi$ between the categories of
algebraic groups and Hopf algebras. This needs the introduction of the new 
notion of {\it crossed comodule} of Hopf algebras (definition \ref{def6}) 
dual to the concept of crossed module of Hopf algebras. The main results of 
this section (propositions \ref{prop3} and \ref{prop4}) together with classical
equivalence between algebraic goups and their rings of functions give a 
correspondence between commutative Hopf $2$-algebras (in their incarnation as 
crossed comodules) and connected algebraic $2$-groups.

\begin{defi}[Definition of a crossed comodule of Hopf algebras.] \label{def6}
Let $\zeta:K\to L$ be a morphism of Hopf algebras. The morphism $\zeta$ defines
a crossed comodule of Hopf algebras in case
\begin{itemize}
\item[(i)] $(L,\rho_L)$ is at the same time a left $K$-comodule algebra and 
a left $K$-comodule coalgebra,
\item[(ii)] $\rho_L\circ\zeta\,=\,(\id_K\otimes\zeta)\circ{\rm coad}_K$,
\item[(iii)] $(\zeta\otimes\id_L)\circ\rho_L\,=\,{\rm coad}_L$.
\end{itemize}
\end{defi}

\begin{defi}
A morphism $(\rho,\sigma):(\zeta:K\to L)\to(\zeta':K'\to L')$ of crossed
comodules of Hopf algebras $\zeta:K\to L$ and $\zeta':K'\to L'$ is a pair
of morphisms of Hopf algebras $\rho:K\to K'$ and $\sigma:L\to L'$ such that 
the following diagrams are commutative:

\vspace{.5cm}
\hspace{3cm}
\xymatrix{
K \ar[r]^{\rho} \ar[d]^{\zeta} & 
K' \ar[d]^{\zeta'} &  & L \ar[r]^{\sigma} \ar[d]^{\rho_L}
& L' \ar[d]^{\rho_{L'}} \\ 
H \ar[r]^{\sigma}  & H' & & K\otimes L \ar[r]^{\rho\otimes\sigma} &
K'\otimes L'  }
\vspace{.5cm}
\end{defi}

It is easy to show that crossed comodules form a category with this notion
of morphisms.

We recall for (i) that, by definition,  $(L,\rho_L)$ is a left $K$-comodule 
coalgebra if $\Delta_L$ is a comodule morphism, i.e. satisfies  $(\id_K\otimes
\Delta_L)\circ\rho_L\overset{\Diamond}{=}\rho_{L\otimes L}\circ\Delta_L$ with 
$\rho_{L\otimes L}:=
(\mu_K\otimes \id_L\otimes \id_L)\circ(\id_K\otimes\tau\otimes \id_L)\circ
(\rho_{L}\otimes \rho_{L}).$

In (ii) and (iii) ${\rm coad}_L$ and ${\rm coad}_K$ denote the adjoint 
coaction (and not the coadjoint action !) defined as ${\rm coad}_L\,:=\,
(\mu_L\otimes\id_L)\circ(\id_L\otimes\tau)\circ(\triangle_L
\otimes S_L)\circ\triangle_L.$
  
The need for definition \ref{def6} relies on the fact that the notion of a 
Hopf algebra is auto-dual, while that of a crossed module is not. 
The existence of both notions of crossed modules and comodules shows that 
there are (at least) two natural ways to
categorify the notion of a Hopf algebra. \\

The next proposition generates natural examples of crossed comodules of Hopf 
algebras, namely when one considers the space of functions on an algebraic 
$2$-group.

\begin{prop}\label{prop3}
The functor $k[-]$ extends to a functor (still denoted $k[-]$) from the 
category of crossed modules of connected algebraic groups to the category of 
crossed comodules of Hopf algebras.
\end{prop}

\begin{proof} Let $\mu:M\to N$ be a crossed module of connected algebraic 
groups. We get a 
morphism of commutative Hopf algebras $\zeta:=\mu^*:k[N]\to k[M]$. 
The action of $N$ on $M$ is transformed into a coaction of $k[N]$ on $k[M]$: 
$\rho_{k[M]}:k[M]\to k[N]\otimes k[M]$ is just the dual of the action map 
$N\times M\to M$. Now, $k[M]$ is a $k[N]$-comodule coalgebra, because of the
fact that $N$ acts on $M$ by automorphisms. On the other hand, $k[M]$ is a 
$k[N]$-comodule algebra, because the multiplication on $k[M]$ is the dual 
of the diagonal $M\to M\times M$, $m\mapsto (m,m)$. 

Finally, the identities (ii) and (iii) of crossed comodule of Hopf 
algebras are consequences (dualization) of the 
properties (a) and (b) of the crossed module $\mu:M\to N$.

The morphism-part of the statement is straight forward.\end{proof} 

We didn't say much about the equivalence of identities (ii) and (iii) and 
properties (a) and (b) here: one point of view is that (a) and (b) can be 
formulated in terms of morphisms (suppressing objects), and then this 
equivalence is a formal duality, i.e. one simply applies the functor $k[-]$. 
Another point of view is to write everything in terms of objects -- this has 
the advantage of being more explicit. We will follow this line in the proof of 
proposition \ref{prop4} below.  

Now, conversely, one can obtain obtain $2$-groups from crossed comodules of 
Hopf algebras. This is the content of the following proposition.

\begin{prop}\label{prop4}
The functor $\chi$ extends to a functor (still denoted $\chi$) from the 
category of crossed comodules of commutative Hopf algebras to the category of
crossed modules of groups.
\end{prop}
Its proof will be given after having introduced some notations and the 
preliminary lemma \ref{lem4}.

For the rest of this section, $K$ and $L$ will be two Hopf algebras,  
$\rho_L:L\mapsto K\otimes L$ a linear map, $M:=\chi[L]$, 
$N:=\chi[K]$, $\eta_K,\gamma_K\in N$, $\phi_L$, $\psi_L\in M$. 

One defines the product $\star_{\rho_L}$ of $N$ on $M$ by 
$$\eta_K\star_{\rho_L}\phi_L:=\mu_k\circ(\eta_K\otimes\phi_L)\circ\rho_L.$$ 
In case $K=L$ and $\rho_L=\Delta_L$ (resp. $L=K$ and $\rho_L=\Delta_K$, 
resp. $L=K$ and $\rho_L={\rm coad}_K$), we will simplify the notation 
$\star_{\rho_L}$ to $\star_L$ (resp. $\star_K$, resp. 
$\star_{\overset{\sim}{\rho_K}}$). From now on, $(L, \rho_L)$ will be a 
left $K$-comodule algebra and coalgebra.

\begin{lem}\label{lem4}
\begin{itemize}
\item[($\alpha$)] $(M, \star_L)$ is a group with inverse $S_L^\star$,
\item[($\beta$)] $\star_{\rho_L}$ defines a group action of $N$ on $M$,
\item[($\gamma$)] $\star_{\rho_L}$ is an action by automorphisms,
\item[($\delta$)]$\eta_K\star_{\overset{\sim}{\rho_K}}\gamma_K=\eta_K\star
\gamma_K\star\eta_K^{-1}$.
\end{itemize}
\end{lem}
\begin{proof} ($\alpha$) and ($\beta$) are proven using definitions and 
coassociativity 
of $\Delta_L$ and $\rho_L$. The proof of ($\gamma$), i.e. 
$N$ acts on $M$ by automorphisms, is a little more involved and is based on 
the fact ($\Diamond$) that $\Delta_L$ is a comodule map:
\begin{eqnarray*}
\eta_K\star_{\rho_L}(\phi_L\star_L\psi_L) & \overset{\star_{\rho_L}}{=}  & 
\mu_k(\eta_K\otimes(\phi_L\star_L\psi_L))\circ\rho_L \\
                                          & \overset{\star_L}{=}  & 
\mu_k(\id_k\otimes\mu_k)\circ(\eta_K\otimes\phi_L\otimes\psi_L)\circ
(\id_K\otimes\Delta_L)\circ\rho_L\\
                                          & \overset{(\Diamond)}{=}  & 
\mu_k(\id_k\otimes\mu_k)\circ(\eta_K\otimes\phi_L\otimes\psi_L)\circ
\rho_{L\otimes L}\circ\Delta_L\\
                                           & \overset{\rho_{L\otimes L}}{=}  & 
\mu_k(\id_k\otimes\mu_k)\circ(\eta_K\otimes\phi_L\otimes\psi_L) \\
& & \circ(\mu_K\otimes \id_L\otimes \id_L)\circ(\id_K\otimes\tau\otimes \id_L)
\circ(\rho_{L}\otimes \rho_{L})\circ\Delta_L\\
                                          & \overset{\eta_K}{=}  & 
\mu_k(\id_k\otimes\mu_k)\circ(\mu_k(\eta_K\otimes\eta_K)\otimes\phi_L\otimes
\psi_L) \\
& & \circ(\id_K\otimes\tau\otimes \id_L)\circ(\rho_{L}\otimes \rho_{L})\circ
\Delta_L\\
                                          & \overset{\tau}{=}  & 
\mu_k(\mu_k\otimes\mu_k)\circ(\eta_K\otimes\phi_L\otimes\eta_K\otimes\psi_L)
\circ(\rho_{L}\otimes \rho_{L})\circ\Delta_L\\
                                            & \overset{\star_L}{=}  & 
(\mu_k\circ (\eta_K\otimes\phi_L)\circ\rho_{L})\star_L(\mu_k\circ (\eta_K
\otimes\psi_L)\circ\rho_{L})\\
                                           & \overset{\star_{\rho_L}}{=}  & 
 (\eta_K\star_{\rho_L}\phi_L)\star_L(\eta_K\star_{\rho_L}\psi_L).\\
\end{eqnarray*}
We finally prove ($\delta$):
\begin{eqnarray*}
\eta_K\star_{\rho_K}\gamma_K & \overset{\star_{\rho_K}}{=} & \mu_k(\eta_K
\otimes\gamma_K)\circ {\rm coad}_K\\
                                  & \overset{coad_K}{=} & \mu_k(\eta_K\otimes
\gamma_K)\circ (\mu_K\otimes\id_K)\circ(\id_K\otimes\tau)\circ(\triangle_K
\otimes S_K)\circ\triangle_K\\
                                & \overset{\mu_K}{=} & \mu_k(\mu_k\otimes\id_k)
\circ(\eta_K\otimes\eta_K\otimes\gamma_K)\circ(\id_K\otimes\tau)\circ(
\triangle_K
\otimes S_K)\circ\triangle_K\\
                                & \overset{\tau}{=} & \mu_k(\mu_k\otimes\id_k)
\circ(\eta_K\otimes\gamma_K\otimes\eta_K)\circ(\triangle_K
\otimes S_K)\circ\triangle_K\\
                                & \overset{\star_K}{=} &\eta_K\star\gamma_K
\star\eta_K^{-1}. \\
\end{eqnarray*}
\end{proof} 
We can now turn to the proof of proposition \ref{prop4}.\\

\begin{proof} Let $\zeta:K\to L$ be a crossed comodule of commutative Hopf 
algebras
$K$ and $L$. 
By functoriality, we get a group homomorphism 
$\mu:=\zeta^*:M\to N$.  We already know from lemma \ref{lem4} that 
$\star_{\rho_L}$ is an action by automorphisms, so it only remains to show 
that identities 
(a) and (b) hold. \\

Identity (a) is obtained from axiom (ii): 
\begin{eqnarray*}
\mu (\eta_K\star_{\rho_L}\psi_L) & \overset{\mu}{=} & (\eta_K\star_{\rho_L}
\psi_L)\circ\zeta \\
                                & \overset{\star_{\rho_L}}{=} & \mu_k(\eta_K
\otimes\psi_L)\circ\rho_L\circ\zeta \\
                                & \overset{(ii)}{=} & \mu_k(\eta_K\otimes
\psi_L)\circ (\id_K\otimes\zeta)\circ {\rm coad}_K \\
                                & \overset{\mu}{=} & \mu_k(\eta_K\otimes\mu(
\psi_L))\circ {\rm coad}_K \\
                                & \overset{\star_{\overset{\sim}{\rho_K}}}{=} 
&(\eta_K\star_{\rho_K}\mu(\psi_L)) \\
\end{eqnarray*}
Identity (b) is obtained from axiom (iii) and lemma \ref{lem4} ($\rho$):
\begin{eqnarray*}
\mu (\phi_L)\star_{\rho_L}\psi_L & \overset{\mu}{=} & (\phi_L\circ\zeta) 
\star_{\rho_L}\psi_L   \\
                                & \overset{\star_{\rho_L}}{=} & \mu_k(\phi_L
\otimes\psi_L)\circ(\zeta\otimes id_L)\circ\rho_L \\
                                & \overset{(iii)}{=} & \mu_k(\phi_L\otimes
\psi_L)\circ {\rm coad}_L \\
                                & \overset{(\delta)}{=} & \phi_L\star\psi_L
\star\phi_L^{-1}. \\
\end{eqnarray*}
\end{proof} 

\section{Integration of Lie $2$-algebras}

In this section, we explain a general scheme of how to integrate (some class 
of) semi-strict Lie $2$-algebras. This is motivated by the following example 
of a crossed module of Hopf algebras:

\begin{rem} \label{rem5}
Let us treat in some detail the example of the crossed module of Hopf algebras 
corresponding to the crossed module of Lie algebras which represents the
generator of $H^3({\mathfrak g},\C)$, where ${\mathfrak g}$ is
any simple complex finite dimensional Lie algebra. This leads
to the definition of the enveloping $2$-algebra of the string Lie algebra
of Baez and Crans \cite{BaeCra}.

Denote by ${\mathfrak h}$ a Cartan subalgebra of ${\mathfrak g}$ and choose a
Borel subalgebra ${\mathfrak b}\supset{\mathfrak h}$. Let ${\mathfrak n}$
be the nilpotent subalgebra of ${\mathfrak g}$ such that ${\mathfrak n}\oplus
{\mathfrak h}={\mathfrak b}$ as vector spaces.
Denote by $M(\lambda)$ the Verma module of ${\mathfrak g}$ of highest weight
$\lambda$. 
Namely $M(\lambda)=U{\mathfrak g}\otimes_{U{\mathfrak b}}\C_{\lambda}$,
where $\C_{\lambda}$ is the one dimensional ${\mathfrak b}$-module given by
the trivial action of ${\mathfrak n}$ and the ${\mathfrak h}$-action via
$\lambda\in{\mathfrak h}^*$. $M(\lambda)$ possesses a unique maximal proper
submodule $N(\lambda)$ and the quotient $L(\lambda)$ is therefore irreducible.
Denote by $M(\lambda)^{\sharp}$, $N(\lambda)^{\sharp}$ and 
$L(\lambda)^{\sharp}$ the restricted duals of these graded 
${\mathfrak g}$-modules.
We have by definition a short exact sequence of ${\mathfrak g}$-modules:
$$0\to L(\lambda)^{\sharp}\to M(\lambda)^{\sharp}\to N(\lambda)^{\sharp}\to 0,
$$
which specializes to 
$$0\to \C\to M(0)^{\sharp}\to N(0)^{\sharp}\to 0$$
for $\lambda=0$. On the other hand, there is a $2$-cocycle $\alpha'\in 
C^2({\mathfrak g},N(0)^{\sharp})$ which gives an abelian extension 
$$0\to N(0)^{\sharp}\to N(0)^{\sharp}\times_{\alpha'}{\mathfrak g}\to 
{\mathfrak g}\to 0.
$$
It is shown in \cite{Wag} that the splicing together of these two sequences 
gives a crossed module
$$\mu:N(0)^{\sharp}\to N(0)^{\sharp}\times_{\alpha'}{\mathfrak g},$$
which represents a generator of $H^3({\mathfrak g};\C)$.

Another crossed module representing a generator of $H^3({\mathfrak g};\C)$ is
known \cite{Nee}: denoting by $P{\mathfrak g}$ and $\Omega{\mathfrak g}$ the 
path- and loop Lie algebra corresponding to ${\mathfrak g}$, there is a 
(general) extensions of Lie algebras 
$$0\to\Omega{\mathfrak g}\to P{\mathfrak g}\to{\mathfrak g}\to 0,$$
adapting the path-loop fibration to the Lie algebra setting. But a crossed 
module is nothing but a central extension of an ideal of some Lie algebra 
(cf \cite{Nee}), and thus the standard central extension of 
$\Omega{\mathfrak g}$ gives rise to a crossed module
$$0\to\C\to\Omega{\mathfrak g}\to P{\mathfrak g}\to{\mathfrak g}\to 0.$$
It is shown in \cite{Nee} that the corresponding class generates 
$H^3({\mathfrak g};\C)$.  

Recall further from \cite{BaeCra} that semi-strict Lie $2$-algebras are 
categorified Lie algebras with a functorial bracket which is strictly
antisymmetric, but satisfies the Jacobi identity only up to a Jacobiator. 
This Jacobiator gives a $3$-cocycle $\theta$ of one of the underlying Lie 
algebras
${\mathfrak g}$ with values in some ${\mathfrak g}$-module $V$. The data
$({\mathfrak g},V,[\theta])$ then completely specifies the equivalence class
of the given semi-strict Lie $2$-algebra. In the end of their paper 
\cite{BaeCra}, Baez and Crans construct a family ${\mathfrak g}_{\hbar}$
of semi-strict Lie 
$2$-algebras whose equivalence classes are given in this sense by the triplets
$({\mathfrak g},\C,[\hbar\theta])$, where here $\theta$ is the Cartan cocycle,
i.e. $\theta=\langle[,],\rangle$ for $\langle,\rangle$ the Killing form of 
the simple complex finite dimensional Lie algebra ${\mathfrak g}$, and 
$\hbar$ is a scalar.   

Recall now from \cite{BCSS} that the Lie $2$-algebras ${\mathfrak g}_{\hbar}$
are linked to the {\it string group} (for $\hbar=\pm 1$). For a given simply 
connected, simple 
Lie group $G$, there is a topological group $\hat{G}$ obtained by killing the
third homotopy group of $G$. This group $\hat{G}$ is, by analogy with the case
of $G={\rm Spin}(n)$, called the string group of $G$. Now the authors
of \cite{BCSS} construct an infinite-dimensional Lie $2$-group, whose 
Lie $2$-algebra is equivalent to ${\mathfrak g}_{\hbar}$ and whose geometric 
realization is (homotopy equivalent to) $\hat{G}$ for $\hbar=\pm 1$.

Direct application of proposition \ref{prop1} gives us a crossed module of 
Hopf algebras

\begin{equation}  \label{**}
\gamma:S(N(0)^{\sharp})\to S(N(0)^{\sharp})\otimes_{\hbar\alpha'}
U{\mathfrak g}
\end{equation}

Here $S$ denotes the symmetric algebra on a 
$k$-vector space. The above mentioned crossed module corresponding to the 
path-loop fibration also gives rise to an enveloping algebra, but observe that
it is not as algebraic as (\ref{**}) and thus less meaningful (what is the
enveloping algebra of a topological Lie algebra ?) from the point 
of view of representation theory (while it has the advantage of being more 
geometric). More precisely, in our crossed module we take enveloping algebras
of {\it at most countably infinite dimensional} vector spaces, while the 
underlying vector spaces in the crossed module stemming from the
path-loop fibration are not countably infinite dimensional.   
\end{rem}

The construction which we cited here as an exemple works for any $3$-cohomology
class of a finite dimensional Lie algebra. Therefore
we deduce the following theorem (cf \cite{BaeCra} for missing notations)
directly from proposition \ref{prop1}.

\begin{theo}  \label{thm4}
For any skeletal semistrict Lie $2$-algebra corresponding to a triple 
$({\mathfrak g},V,[\theta]\in H^3({\mathfrak g},V))$, there is a 
crossed module of Hopf algebras
\begin{equation}  
\gamma:S(I)\to S(Q)\otimes_{\alpha''}U{\mathfrak g}
\end{equation}
which is a natural {\it algebraic} candidate for the enveloping algebra of the 
given Lie $2$-algebra in the sense that all its underlying vector spaces are
at most countably infinite dimensional.    
\end{theo}

Here $Q$ and $I$ are ${\mathfrak g}$-modules such that $I$ is injective and 
such that there is a short exact sequence of ${\mathfrak g}$-modules:
$$0\to V\to I\to Q\to 0.$$
The class of the $2$-cocycle $\alpha''\in Z^2({\mathfrak g},Q)$ corresponds to 
the class of the given $3$-cocycle $\theta$ under the corresponding 
connecting homomorphism. The above construction (\ref{**}) is a special case.
It is well known that in the context of remark \ref{rem5},
there cannot be a purely finite dimensional construction in general 
(see \cite{Hoch}). 

The properties of the preceeding theorem are important, because they permit to 
establish a duality theory for the mentioned kind of enveloping algebras for 
skeletal semi-strict Lie $2$-algebras (up to equivalence of Lie $2$-algebras,
cf \cite{BaeCra}).

Recall how finite dimensional Lie algebras are linked to formal groups (in
characteristic zero), see \cite{Ser}. Take a (finite dimensional real)
Lie group $G$ with unit $1$, and
complete the algebra of germs of functions around $1$ by the adic topology. 
Its continuous dual may be identified to the Hopf algebra $U$ of 
point distributions.
$U$ is isomorphic to $U{\mathfrak g}$, where ${\mathfrak g}$ is the Lie algebra
of $G$. Synthetically, this can be understood as the conjunction of linear
duality and the Milnor-Moore theorem. In this way, Serre establishes
an equivalence of categories between the category of formal groups and the 
category of finite dimensional Lie algebras.   

Now this also works backwards: starting with a finite dimensional Lie algebra
${\mathfrak g}$, one forms the universal enveloping algebra $U{\mathfrak g}$.
Its dual is the completion of a commutative Hopf algebra, from which one may 
extract the group by taking characters.  

This scheme of reasoning obviously also works for crossed modules of 
Lie algebras; this is our first procedure to integrate Lie $2$-algebras. 
One has to apply the functor $U$ in order to obtain a 
crossed modules of Hopf algebras by proposition \ref{prop1}. Then one passes 
to the continuous dual, which is naturally a crossed comodule of commutative 
Hopf algebras (see section $5$). Each constituant of this crossed comodule
is the completion of the commutative Hopf algebra of functions on some 
formal group, or of germs of functions in an identity-neighborhood of some
Lie group $G$. The results of section $5$ now imply that this crossed 
comodule gives rise to a crossed module of groups, 
which may then be interpreted as some 
$2$-group. In this way, the above theorem permits to integrate 
strict Lie $2$-algebras, i.e. crossed modules of Lie algebras, where all 
constituants are finite dimensional Lie algebras. A general semi-strict
Lie $2$-algebra is always equivalent to a strict one (``strictification''), 
thus this scheme permits 
to integrate semi-strict Lie $2$-algebras which admit a strictification 
consisting of finite dimensional Lie algebras. We summarize the above 
discussion in the following theorem:

\begin{theo}
Every semi-strict Lie $2$-algebra which admits a strictification 
consisting of finite dimensional Lie algebras can be integrated into a 
finite dimensional strict Lie $2$-group whose associated  Lie $2$-algebra
is equivalent to the given one.
\end{theo}

Let us remark that Behrang Noohi has a different point of view which also leads
to an integration procedure for this kind of Lie $2$-algebras, see \cite{Noo}.
  
Let us finally discuss the integrability of the crossed module of Lie algebras
$$\mu:N(0)^{\sharp}\to N(0)^{\sharp}\times_{\alpha'}{\mathfrak g},$$  
to which the theorem does not apply. An integration of the more geometric 
crossed module $\Omega{\mathfrak g}\to P{\mathfrak g}$ is the object of 
\cite{BCSS}. In fact, the crossed module $\mu$ is not integrable (at least, not
in a na\"{\i}ve way).

Indeed, as explained in \cite{Wag}, the crossed module $\mu$ which also 
represents the Cartan class in $H^3{\mathfrak g},\C)$, is spliced together 
from an abelian extension
$$0\to (U{\mathfrak g}^+)^*\to (U{\mathfrak g}^+)^*\times_{\alpha}{\mathfrak g}
\to{\mathfrak g}\to 0,$$
and the short exact sequence of (adjoint) ${\mathfrak g}$-modules
$$0\to\C\to(U{\mathfrak g})^*\to(U{\mathfrak g}^+)^*\to 0.$$ 

Let $G$ be the connected, simply connected compact Lie group corresponding to 
${\mathfrak g}$. The non-integrability of $\mu$ is due to the fact that
the ${\mathfrak g}$-module $U{\mathfrak g}$ (with the action 
by left multiplication) does not integrate into a $G$-module. 

Even if it did, the integration would not have been this simple: in fact,
the resulting crossed module of groups must contain a topologically non-trivial
fiber bundle. If not, there would be a globally smooth section of the crossed
module, and its $3$-class would be presented by a globally smooth cocycle. 
This is impossible, because the smooth cohomology of $G$ is trivial (see
for example the references to Hu's, van Est's or Mostow's work in 
\cite{BaeLau}).

Notwithstanding this example, in principle crossed modules of the above kind 
which are spliced together from an abelian extension and a short exact sequence
of modules may be integrated as follows -- this is our (sketched) second 
integration procedure, which applies (in principle) to infinite dimensional Lie
algebras: first integrate the modules to 
locally convex $G$-modules, and suppose that the $2$-cocycle from the abelian 
extension is continuous. Integrate it into a locally smooth group $2$-cocycle 
using \cite{Neeb}. Obstructions to this step reside in $\pi_1(G)$ and 
$\pi_2(G)$ (see {\it loc. cit.}) and vanish therefore. Finally, splice together
the exact sequences to obtain a crossed module of groups.

\section{Relation to previous work on crossed modules of Hopf algebras}

The only work on crossed modules of Hopf algebras we know about is
\cite{FerLopNov}. Its main result, theorem (14), states the 
equivalence of categories between the category of ${\rm cat}^1$--Hopf algebras 
and the category of crossed modules of Hopf algebras. 
The definition of pre-${\rm cat}^1$--Hopf 
algebras that it offers can be generalised as follows:

\begin{defi} \label{def7}
Let $A$ and $H$ be two Hopf algebras, 
$s,t:A\to H$ and $e:H\to A$ be morphisms of 
Hopf algebras. $\overline{A}:=(A,H,s,t,e)$ is a pre-${\rm cat}^1$-Hopf 
algebra if $$s\circ e= t\circ e=id_H.$$
\end{defi}

\begin{rem} \label{rem10}
The definition of a pre-${\rm cat}^1$-Hopf 
algebra given in \cite{FerLopNov} requires in addition cocommutativity of $H$ 
together with the conditions ii) $(s\otimes A)\circ \tau_{A,A}\circ \Delta_A=
(s\otimes A)\circ\Delta_A$ and iii)$(t\otimes A)\circ  \tau_{A,A}\circ 
\Delta_A=(t\otimes A)\circ\Delta_A$.
\end{rem}

Then, in definition (9) of \cite{FerLopNov}, ${\rm cat}^1$-Hopf algebras are 
defined as monoids in 
a certain monoidal category $\mathcal{P}C^1_H$ whose objects are precisely  
pre-${\rm cat}^1$-Hopf 
algebras in the sense of \cite{FerLopNov} (we will not recall the monoidal 
structure here).

The proof of the equivalence between the category of ${\rm cat}^1$-Hopf 
algebras and the category of crossed modules of Hopf algebras (theorem (14) 
in \cite{FerLopNov}) is rather 
technical and involves abstract tools taken from other works. We propose here 
a different approach to this equivalence.

\begin{theo}  \label{thm6}
When restricting to the subcategory of irreducible cocommutative Hopf algebras,
the categories of crossed modules of Hopf algebras and the category of 
pre-${\rm cat}^1$-Hopf algebra are equivalent.
\end{theo}

\begin{proof} 
Given a pre-${\rm cat}^1$-Hopf algebra $s,t:A\to H$, $e:H\to A$ where $H$ 
and $A$ are 
irreducible cocommutative Hopf algebras (over the characteristic zero field 
$k$), we have $A=U(P(A))$ and $H=U(P(H))$ (see \cite{Swe} theorem 13.0.1, 
p. 274), and the morphisms are also induced from morphisms of Lie algebras.
Denote by ${\mathfrak g}_1:=P(A)$ and ${\mathfrak g}_0:=P(H)$, then we have a 
Lie $2$-algebra ${\mathfrak g}_1\to{\mathfrak g}_0$, because the 
categorical composition $m$ is already encoded in the $2$-vector space 
structure (see section $3$). 
This Lie $2$-algebra corresponds then to a crossed module of Lie 
algebras $\mu:{\mathfrak m}\to{\mathfrak n}$, and applying the functor $U$,
we get the crossed module of Hopf algebras associated to the given 
pre-${\rm cat}^1$-Hopf algebra $s,t:A\to H$ by proposition $1$.

Conversely, given a crossed module of irreducible cocommutative Hopf algebras
$\gamma:B\to H$, we can associate to it in the same way a 
pre-${\rm cat}^1$-Hopf algebra $s,t:A\to H$, $e:H\to A$ using proposition $2$.
\end{proof}

\begin{rem} \label{rem11}
In the same way, one can show that the category of crossed comodules of 
commutative Hopf algebras and the category of commutative Hopf $2$-algebras 
(meaning explicitly pre-${\rm cat}^1$-Hopf algebras with commutativity 
replacing the cocommutativity requirement) are equivalent when restricting 
to the subcategory of those Hopf algebras which are finitely generated as an 
algebra. It suffices to apply the functors $k[-]$ and $\chi$ to translate the 
problem into groups, where we then apply theorem \ref{thm2}.
\end{rem}
  
\begin{rem}  \label{rem12}
In conjunction with theorem $(14)$ of {\it loc. cit.}, our proposition shows 
that when one restricts the monoidal category $\mathcal{P}C^1_H$ to 
irreducible cocommutative Hopf algebras, every object becomes a monoid. 
In particular, this implies that the characterization of monoids given in 
proposition (\cite{FerLopNov}, 14 ii)) should be automatically satisfied.
\end{rem}

\begin{rem}  \label{rem13}
The previous remark can be seen as the incarnation in the Hopf level of 
remark \ref{rem4}. A natural question at this stage is whether this property 
holds true if one drops the irreducibility property. In the terms of 
\cite{FerLopNov}, is every object of $\mathcal{P}C^1_H$ a monoid ?
\end{rem}

\end{document}